\documentclass[11pt,a4paper]{article}
\usepackage{soul}
\usepackage{amssymb}
\usepackage{dsfont}
\usepackage{graphicx} 
\usepackage{amsmath}
\usepackage{amsthm}
\usepackage{enumitem}   
\usepackage[dvipsnames]{xcolor}
\usepackage{hyperref}
\hypersetup{
    colorlinks=true,
    linkcolor=Green,
    citecolor=Green,
    urlcolor=Cyan
    }
\usepackage{soul}
\usepackage{diagbox}
\usepackage{multirow}
\usepackage{subcaption}
\usepackage[titletoc,title]{appendix}
\usepackage{titlesec}
\usepackage{pgfplotstable}
\usepackage{pgfplots}
\usepgfplotslibrary{fillbetween}
\pgfplotsset{compat = newest}

\titleformat{\paragraph}
{\normalfont\normalsize\bfseries}{\theparagraph}{1em}{}
\titlespacing*{\paragraph}
{0pt}{3.25ex plus 1ex minus .2ex}{1.5ex plus .2ex}

\def\num#1{\numx#1}\def\numx#1e#2{{#1}\mathrm{e}{#2}}
\newcommand*{\collauthor}[2]{{#1}$^{#2}$}

\newcommand*{\affiliation}[2]{$\mbox{}^{{#2}}${#1}}
\newcommand*{\colltitle}[1]{\textbf{#1}}

\newenvironment{keywords}[1]{\vspace{1cm}\\{\bf \slshape{Keywords}}\quad\slshape{#1}}{}
\DeclareMathOperator*{\argmin}{arg\,min}
\newtheorem{example}{Example}[]

\newtheorem{remark}{Remark}[section]

\DeclareMathOperator{\Tr}{Tr}
\DeclareMathOperator{\Hessian}{Hess}

\newcommand{\Kf}{\mathfrak{K}}

\newcommand{\qfrac}{\mathfrak{q}}

\def\num#1{\numx#1}\def\numx#1e#2{{#1}\mathrm{e}{#2}}

\begin{document}
\begin{center}
\begin{Large}
  \colltitle{A forward differential deep learning-based algorithm for solving high-dimensional nonlinear backward stochastic differential equations}
\end{Large}
\vspace*{1.5ex}

\begin{sc}
\begin{large}
\collauthor{Lorenc Kapllani}{}
and
\collauthor{Long Teng}{}
\end{large}
\end{sc}
\vspace{1.5ex}

\affiliation{Chair of Applied and Computational Mathematics,\\
Faculty of Mathematics and Natural Sciences,\\
University of Wuppertal,\\
Gau{\ss}str. 20, 42119 Wuppertal, Germany\linebreak }{} \\
\end{center}

\section*{Abstract}
In this work, we present a novel forward differential deep learning-based algorithm for solving high-dimensional nonlinear backward stochastic differential equations (BSDEs). Motivated by the fact that differential deep learning can efficiently approximate the labels and their derivatives with respect to inputs, we transform the BSDE problem into a differential deep learning problem. This is done by leveraging Malliavin calculus, resulting in a system of BSDEs. The unknown solution of the BSDE system is a triple of processes $(Y, Z, \Gamma)$, representing the solution, its gradient, and the Hessian matrix. The main idea of our algorithm is to discretize the integrals using the Euler-Maruyama method and approximate the unknown discrete solution triple using three deep neural networks. The parameters of these networks are then optimized by globally minimizing a differential learning loss function, which is novelty defined as a weighted sum of the dynamics of the discretized system of BSDEs. Through various high-dimensional examples, we demonstrate that our proposed scheme is more efficient in terms of accuracy and computation time compared to other contemporary forward deep learning-based methodologies.
\begin{keywords}
backward stochastic differential equations, high-dimensional problems, deep neural networks, differential deep learning, global optimization, local loss function, malliavin calculus, option pricing and hedging
\end{keywords}

\section{Introduction}
\label{sec1}
In this work, we develop a novel forward differential deep learning-based scheme to numerically solve decoupled high-dimensional forward-backward stochastic differential equations (FBSDEs)
\begin{equation}
    \begin{split}
        \left\{
            \begin{array}{rcl}
                X_t & = & x_0 + \int_{0}^{t} a \left(s, X_s\right)\,ds + \int_{0}^{t} b \left(s, X_s\right)\,dW_s,\\
   	   		    Y_t & = & g\left(X_T\right) + \int_{t}^{T} f\left(s, \mathbf{X}_s\right)\,ds -\int_{t}^{T}Z_s\,dW_s,
            \end{array} \quad \forall\, t \in [0, T]
        \right. 
    \end{split}
\label{eq1}
\end{equation}
where $\mathbf{X}_t := \left( X_t, Y_t, Z_t\right)$, $W_t = \left( W_t^1, \ldots, W_t^d \right)^\top$ is a $d$-dimensional Brownian motion, $a: [0, T] \times \mathbb{R}^{d} \to \mathbb{R}^{d}$, $b: [0, T] \times \mathbb{R}^{d} \to \mathbb{R}^{d \times d}$, $f:\left[0,T\right]\times\mathbb{R}^{d}\times\mathbb{R}\times\mathbb{R}^{1\times d} \to\mathbb{R}$ is the driver function and $g: \mathbb{R}^{d} \to \mathbb{R}$ is the terminal condition. The terminal condition $Y_T$ depends on the final value $X_T$ of the forward stochastic differential equation (SDE). Usually, the coupled FBSDE is referred to as a FBSDE. Hence, to avoid confusion, we refer to the decoupled FBSDE~\eqref{eq1} as a BSDE in the remainder of this paper.

The existence and uniqueness of the solution of~\eqref{eq1} is proven in~\cite{Pardoux1990}. After that, BSDEs have found many applications across various scientific domains, such as finance and physics, due to their connection to partial differential equations (PDEs) through the nonlinear Feynman-Kac formula. In finance, the solution $(Y, Z)$ of a BSDE provides the price and delta-hedging, see~\cite{El1997}. In many practical applications, BSDEs are often nonlinear and high-dimensional, and analytical solutions are typically not available. Hence, advanced numerical techniques to approximate their solution become interesting. In the recent years, many numerical methods have been proposed for solving BSDEs, we refer to~\cite{chessari2023numerical} for a nice overview of the classical approaches, e.g. Fourier or cubature methods on spatial discretization. However, most of them suffer from the ``curse of dimensionality", where the computational cost increases exponentially with the problem's dimensionality. 

To address this challenge, several works have introduced innovative algorithms for solving high-dimensional nonlinear BSDEs, which can be classified into three main categories. The first category involves multilevel Monte Carlo methods based on Picard iteration~\cite{weinan2019multilevel,becker2020numerical,hutzenthaler2020overcoming,hutzenthaler2020multilevel,hutzenthaler2021multilevel,hutzenthaler2022overcoming,hutzenthaler2023overcoming}. The second category includes tree-based methods~\cite{crisan2012solving,teng2021review,teng2022gradient}, and the third one consists of deep learning-based methods using deep neural networks (DNNs)~\cite{weinan2017deep,han2018solving,fujii2019asymptotic,hure2020deep,takahashi2022new,beck2021deep,germain2022approximation,andersson2023convergence,gnoatto2023deep,raissi2024forward,kapllani2024deep}. Recently, a novel category of schemes has been proposed, referred to as differential deep learning-based schemes~\cite{kapllani2024backward}, which can be considered as a generalization of deep learning-based methods. The latter has been shown to outperform the deep learning scheme~\cite{hure2020deep} in approximating the process $Y$ and especially the processes $( Z,\Gamma)$. The triple of processes $\left(Y, Z, \Gamma\right)$ in a BSDE represents the unknown solution, its gradient, and the Hessian matrix. Note that the differential deep learning scheme~\cite{kapllani2024backward} formulates the BSDE as a local optimization problem. However, a differential deep learning scheme based on global optimization that provides high-accurate approximations of $\left(Z, \Gamma\right)$  is missing in the literature. This study aims to fill this gap.

The deep learning schemes are inspired by the pioneering work~\cite{weinan2017deep,han2018solving} called the deep BSDE (we refer to it as the DBSDE scheme). In this approach, the BSDE is formulated as a global optimization problem. After the time domain is discretized, the discrete process $Z$ is parameterized using DNNs. The parameters of DNNs are optimized using the stochastic gradient descent (SGD) algorithm on a loss function defined at the terminal time $T$. A strong drawback of the DBSDE scheme is that it achieves much better approximations of the BSDE~\eqref{eq1} at the initial time than at the other time points, although the solution of the BSDE is approximated pathwise along $[0,T]$. To overcome this drawback, the authors in~\cite{raissi2024forward} introduced a new approach, where the problem is formulated as a global optimization with local loss functions. The process $Y$ is parameterized using a DNN, and its gradient (the process $Z$) is obtained from automatic differentiation (AD). The parameters of the DNN are optimized from the global minimization of the local loss functions defined at each time point, with the loss at terminal time included as an additional term in the loss function. Hence, the proposed algorithm attempt to match the discretized dynamics of the BSDE at each time point. Such schemes that rely on global optimization operate forward in time, we refer to this class as forward deep learning schemes. In contrast, other existing deep learning schemes (e.g.~\cite{hure2020deep}) are based on local optimization, and operate backward in time, which are refereed to as backward deep learning schemes.

Both the forward and backward deep learning schemes often struggle to provide highly accurate first- and second-order gradient approximations, see e.g.~\cite{negyesi2024one,kapllani2024backward} for the reasons. This is crucial for financial applications, particularly in delta- and $\Gamma$-hedging strategies for option contracts. Our work in~\cite{kapllani2024backward} (see also~\cite{negyesi2024one}) uses differential deep learning to improve~\cite{hure2020deep} in the class of backward schemes. Both theoretically and numerically, we demonstrated that our scheme is more efficient compared to~\cite{hure2020deep}, especially in the computation of the processes $\left(Z, \Gamma\right)$. To the best of our knowledge, a forward algorithm to provide high-accurate approximations of $\left(Z, \Gamma\right)$ is missing in the literature. Hence, in this work, we study a novel forward differential deep learning scheme. Unlike~\cite{kapllani2024backward}, our new approach is based on global optimization rather than local optimization.

Our method works as follows. Firstly, we formulate the BSDE as a differential deep learning problem by using Malliavin calculus. This results in a BSDE system, as the Malliavin derivatives of the solution pair $(Y,Z)$ of the BSDE satisfy themselves another BSDE. This formulation requires the estimation of the triple of the processes $(Y, Z, \Gamma)$. Afterward, the BSDE system is discretized using the Euler-Maryuama method and three DNNs are used to parameterize the unknown triple of processes. The parameters of the DNNs are then estimated by globally minimizing a differential learning type loss function, which is novelty defined as a weighted sum of the dynamics of the discretized BSDE system. Hence, the SGD is equipped with explicit information about the dynamics of the processes $\left( Y, Z\right)$ at each discrete time point. As a result, our method can yield more accurate approximations than the forward deep learning counterpart~\cite{raissi2024forward} not only for the process $Z$, but also for the processes $Y$ and $\Gamma$. This is demonstrated in the numerical experiments. Additionally, our algorithm gives significantly shorter computation times compared to~\cite{raissi2024forward} when computing $\Gamma$ at each optimization step. This efficiency is due to the fact that while the latter relies solely on AD to estimate the process $\Gamma$, our method offers the option of using either a DNN or AD, where the use of a DNN proves more time-efficient.

The outline of the paper is organized as follows. In the next section, we introduce some preliminaries including the DNNs. The forward deep learning scheme~\cite{raissi2024forward} and our scheme are presented in Section~\ref{sec3}. Section~\ref{sec4} is devoted to the numerical experiments. Finally, Section~\ref{sec5} concludes this work.

\section{Preliminaries}
\label{sec2}
Let $\left(\Omega,\mathcal{F},\mathbb{P},\{\mathcal{F}_t\}_{0\le t \le T}\right)$ be a complete, filtered probability space. In this space a standard $d$-dimensional Brownian motion $\{W_t\}_{\leq t \leq T}$ is defined, such that the filtration $\{\mathcal{F}_t\}_{0\le t\le T}$ is the natural filtration of $W_t.$ We denote by $| x |$ for the Frobenius norm of any $x \in \mathbb{R}^{d \times \qfrac}$. In the case of scalar and vector inputs, these coincide with the standard Euclidian norm. In what follows, all equalities concerning $\mathcal{F}_t$-measurable random variables are meant in the $\mathbb{P}$-a.s. sense and all expectations (unless otherwise stated) are meant under $\mathbb{P}$. The solution triple $\left(X_t,Y_t,Z_t\right):\left[0,T\right]\times\Omega\to\mathbb{R}^d \times \mathbb{R}\times\mathbb{R}^{1 \times d}$ is the solution of BSDE~\eqref{eq1} if it is $\{\mathcal{F}_t\}_{0\leq t\leq T}$-adapted, square integrable, and satisfies~\eqref{eq1} $\mathbb{P}$-a.s..

An important property of BSDEs is that they provide a probabilistic representation for the solution of a specific class of PDEs given by the nonlinear Feynman-Kac formula. Consider the semi-linear parabolic PDE
\begin{equation}
    \begin{aligned}        
    \frac{\partial u(t,x)}{\partial t} + \nabla_x u(t, x)\,a(t, x) + \frac{1}{2} \Tr\left[b b^{\top} \Hessian_x u (t, x)\right] + f\left(t, x, u, \nabla_x u \,  b\right)(t, x) = 0,
    \end{aligned}
    \label{eq2}
\end{equation}
for all $(t, x) \in ([0, T]\times \mathbb{R}^d)$, the terminal condition $u(T,x)=g(x)$, where $\Hessian_x u$ and $\nabla_x u$ are the Hessian matrix and gradient of function $u$ with respect to spatial variable $x$. Assume that~\eqref{eq2} has a classical solution $u \in \mathbb{C}^{1,2}\left([0, T] \times \mathbb{R}^{d}; \mathbb{R}\right)$ and the regularity conditions of~\eqref{eq1} are satisfied. 
Then the solution of~\eqref{eq1} can be represented $\mathbb{P}$-a.s. by 
\begin{equation}
	Y_t = u\left(t, X_t\right), \quad Z_t= \nabla_x u\left(t, X_t\right)b\left(t, X_t\right) \quad \forall \, t \in \left[0,T\right).
\label{eq3}
\end{equation}
To approximate the function $u$, DNNs are considered due to the approximation capability in high dimensions.

A DNN is a function $\phi(\cdot; \theta): \mathbb{R}^{d_0} \to \mathbb{R}^{d_1}$ composed of a sequence of simple functions, which therefore can be collected in the following form
\begin{equation*}
    x \in \mathbb{R}^{d_0} \longmapsto A_{L+1}(\cdot;\theta(L+1)) \circ \varrho \circ A_{L}(\cdot;\theta(L)) \circ \varrho \circ \ldots \circ \varrho \circ A_1(x;\theta(1)) \in \mathbb{R}^{d_1},
\end{equation*}
where $x \in \mathbb{R}^{d_0}$ is called an input vector, $d_0, d_1\in \mathbb{N}$ is the input and output dimensions, respectively. Moreover $\theta:=\left( \theta(1), \ldots, \theta(L+1) \right) \in \mathbb{R}^{P}$, $P$ is the total number of network parameters and $L \in \mathbb{N}$  is the number of hidden layers each with $\eta \in \mathbb{N}$ neurons. The functions $A_l(\cdot; \theta(l)), l = 1, 2, \ldots, L+1$ are affine transformations: $A_1(\cdot;\theta(1)): \mathbb{R}^{d_0} \to \mathbb{R}^{\eta}$, $A_l(\cdot;\theta(l)),  l = 2, \ldots, L: \mathbb{R}^{\eta} \to \mathbb{R}^{\eta}$ and $A_{L+1}(\cdot;\theta(L+1)): \mathbb{R}^{\eta} \to \mathbb{R}^{d_1}$, represented by
\begin{equation*}
    A_l(v;\theta(l)) = \mathcal{W}_l v + \mathcal{B}_l, \quad v \in \mathbb{R}^{\eta_{l-1}},
\end{equation*}
where $\theta(l):=\left(\mathcal{W}_l, \mathcal{B}_l\right)$, $\mathcal{W}_l \in \mathbb{R}^{\eta_{l} \times \eta_{l-1}}$ is the weight matrix and $\mathcal{B}_l \in \mathbb{R}^{\eta_{l}}$ is the bias vector with $\eta_0 = d_0, \eta_{L+1} = d_1, \eta_l = \eta$ for $l = 1, \ldots, L$. Finally, $\varrho: \mathbb{R} \to \mathbb{R}$ is a nonlinear function (called the activation function), and applied component-wise on the outputs of $A_l(\cdot;\theta(l))$. Common choices are $\tanh(\cdot), \sin(\cdot), \max(0,\cdot)$ etc. We denote by $\Theta$ the set of possible parameters for the DNN $\phi(\cdot; \theta)$ with $\theta \in \Theta$. The universal approximation theorem~\cite{hornik1989multilayer,cybenko1989approximation} justifies the use of DNNs as function approximators.

\section{A forward differential deep learning scheme for BSDE}
\label{sec3}
In this section, we review the forward deep learning scheme~\cite{raissi2024forward} and introduce our new scheme based on differential deep learning. 
\subsection{The local deep BSDE scheme}
\label{subsec31}
The authors in~\cite{raissi2024forward} proposed to formulate the BSDE problem based on a global optimization with
local losses (we refer as Local DBSDE in the rest of the paper). 

The first step is to discretize the integrals in the BSDE~\eqref{eq1}. 
Let us consider $\Delta = \{t_0, t_1, \ldots, t_N\}$ as the time discretization of $[0, T]$ with $t_0 = 0 < t_1 < \ldots < t_N = T$, $\Delta t = t_{n+1} - t_n$. For notational convenience we write $\Delta W_n = W_{t_{n+1}} - W_{t_n}$, $(X_n, Y_n, Z_n) = (X_{t_n}, Y_{t_n}, Z_{t_n})$ and $\left(X^{\Delta}_n, Y^{\Delta}_n, Z^{\Delta}_n\right)$ for the approximations. Applying the Euler-Maruyama scheme in~\eqref{eq1} yeilds
\begin{equation}
     X^{\Delta}_{n+1} = X^{\Delta}_n + a\left(t_n, X^{\Delta}_n\right) \Delta t + b\left(t_n, X^{\Delta}_n\right) \Delta W_n,
     \label{eq4}
\end{equation}
for $n = 0, 1, \ldots, N-1$, $X^{\Delta}_0 = x_0$, and
\begin{equation}
     Y^{\Delta}_n = Y^{\Delta}_{n+1} + f\left(t_n, \mathbf{X}_n^{\Delta}\right) \Delta t -  Z^{\Delta}_n \Delta W_n,
     \label{eq5}
\end{equation}
for $n = N-1, N-2, \ldots, 0,$ where $\mathbf{X}_n^{\Delta} := \left( X_n^{\Delta}, Y_n^{\Delta}, Z_n^{\Delta}\right)$ and $Y^{\Delta}_N = g\left(X^{\Delta}_N\right)$. 

After discretizing the integrals, the scheme in~\cite{raissi2024forward} is made fully implementable by approximating the unknown processes $(Y_n^{\Delta}, Z_n^{\Delta})$ in~\eqref{eq5} for $n=0, 1, \ldots, N$. Due to~\eqref{eq3}, a DNN is used to approximate $Y_n^{\Delta}$ and AD for $Z_n^{\Delta}$. More precisely, the LDBSDE scheme works as follows:
\begin{itemize}
    \item Generate approximations $X^{\Delta}_{n+1}$ for $n = 0, 1, \ldots, N-1$ using~\eqref{eq4}.
    \item At each discrete time point $t_n$, $n = 0, 1, \ldots, N$, use DNN $\phi^y(\cdot; \theta): \mathbb{R}^{1+d} \to \mathbb{R}$ to approximate $Y_n^{\Delta}$ and $Z_n^{\Delta}$ using AD due to \eqref{eq3}, where the input vector of the network is the time value $t_n \in \mathbb{R}_{+}$ and the Markovian process $X_n^{\Delta} \in \mathbb{R}^d$, namely
    $$Y_n^{\Delta,\theta} = \phi^y(t_n, X_n^{\Delta};\theta), \quad Z_n^{\Delta,\theta} = \nabla_x \phi^y(t, x;\theta)\Bigr|_{(t, x) = (t_n,X_n^{\Delta})} b \left( t_n, X_n^{\Delta} \right).$$
    \item Train the parameters $\theta$ using a global loss function including local losses such that the dynamics of discretized BSDE~\eqref{eq5} are satisfied at each time step. The loss is given as    \begin{equation}
        \mathbf{L}^{y,\Delta}\left( \theta \right):= \mathbb{E}\left[\sum_{n=0}^{N-1} \left\vert Y^{\Delta, \theta}_{n+1} - Y_n^{\Delta, \theta} + f\left(t_n, \mathbf{X}^{\Delta, \theta}_n\right) \Delta t -  Z_n^{\Delta,\theta} \Delta W_n  \right\vert^2 + \left\vert Y_N^{\Delta, \theta} - g\left(X_N^{\Delta}\right)\right\vert^2\right],
        \label{eq6}
    \end{equation}
     where for notational convinience $\mathbf{X}^{\Delta, \theta}_n:=\left( X_n^{\Delta}, Y_n^{\Delta, \theta}, Z_n^{\Delta, \theta} \right)$.
     \item Approximate the optimal parameters $\theta^* \in \argmin_{\theta \in \Theta} \mathbf{L}^{y,\Delta}\left( \theta \right)$ using a SGD method and receive the final estimated parameters $\hat{\theta}$. Set the final approximation of $\left( Y_n^{\Delta}, Z_n^{\Delta} \right)$ as $\left( Y_n^{\Delta, \hat{\theta}}, Z_n^{\Delta, \hat{\theta}}\right)$ for $n=0,1,\ldots,N$. 
\end{itemize}

\subsection{The differential local deep BSDE scheme}
\label{subsec32}
To improve the accuracy of first- and second-order gradient approximations -- namely the processes $Z$ and $\Gamma$ -- in the LDBSDE scheme, we use differential deep learning~\cite{huge2020differential}. This limitation in the LDBSDE scheme becomes apparent in its loss function~\eqref{eq6}, as the SGD algorithm lacks the explicit information about the dynamics of $Z$ and does not explicitly include $\Gamma$. In a differential deep learning problem, the loss function requires explicit information about the labels and their derivatives with respect to inputs. Therefore, transforming the BSDE into a differential deep learning problem provides the necessary information to the SGD algorithm. This is done by using the Malliavin calculus.

Applying the Malliavin derivative to~\eqref{eq1} yields another BSDE given as (see~\cite{kapllani2024backward})
\begin{equation}
    \begin{split}
        \left\{
            \begin{array}{rcl}
                D_s X_t & = & \mathds{1}_{s \leq t}\Bigl[ b\left(s, X_s\right)  + \int_s^t \nabla_x a\left( r, X_r \right) D_s X_r dr +  \int_s^t \nabla_x b\left( r, X_r \right) D_s X_r dW_r\Bigr],\\
   	   	 D_s Y_t & = & \mathds{1}_{s \leq t}\Bigl[\nabla_x g\left(X_T\right) D_s X_T+ \int_t^T f_D\left( r, \mathbf{X}_r, \mathbf{D}_s \mathbf{X}_r\right) dr -  \int_t^T \left(\left(D_s Z_r\right)^{\top} dW_r\right)^{\top} \Bigl],
         \end{array}
        \right. 
    \end{split}
\label{eq7}
\end{equation}
where we introduced the notations $\mathbf{D}_s\mathbf{X}_t := \left( D_s X_t, D_s Y_t, D_s Z_t\right)$ and $f_D\left(t, \mathbf{X}_t, \mathbf{D}_s\mathbf{X}_t \right):= \nabla_x f\left( t, \mathbf{X}_t \right) D_s X_t + \nabla_y f\left( t, \mathbf{X}_t\right) D_s Y_t + \nabla_z f\left( t, \mathbf{X}_t\right) D_s Z_t$ $\forall \, 0 \leq s, t \leq T$. Note that $\left(D_s X_t, D_s Y_t, D_s Z_t\right)$ represents the Maliavin derivative of $\left( X_t, Y_t, Z_t\right)$ at time $s \in \left[0, T\right]$. Moreover, $D_t Y_t$ defined by the above equation is a version of $Z_t$, i.e.,
\begin{equation}
    D_t Y_t = Z_t,
    \label{eq8}
\end{equation}
$\mathbb{P}$-a.s. $\forall \, t \in [0, T]$, see~\cite{El1997}. The solution to BSDE system~\eqref{eq1} and~\eqref{eq7} is a pair of triples of stochastic processes $\{\left(X_t, Y_t, Z_t\right)\}_{0\leq t \leq T}$ and $\{\left(D_s X_t, D_s Y_t, D_s Z_t\right)\}_{0\leq s, t \leq T}$ such that~\eqref{eq1} and~\eqref{eq7} holds $\mathbb{P}$-a.s.

As for the LDBSDE scheme, we firstly discretize the integrals in BSDE system~\eqref{eq1} and~\eqref{eq7}. For BSDE~\eqref{eq1}, this is given in~\eqref{eq4} and~\eqref{eq5}. The Euler-Maruyama method gives the approximations of the Malliavin derivatives $D_sX_t$ and $D_s Y_t$ in~\eqref{eq7} as
\begin{equation}
    D_n X_m^{\Delta} = 
    \begin{cases}
        \mathds{1}_{n = m}b\left( t_n, X_n^{\Delta}\right), \quad  0 \leq m \leq n \leq N,\\
   	D_n X^{\Delta}_{m-1} + \nabla_x a\left(t_{m-1}, X^{\Delta}_{m-1}\right) D_n X_{m-1}^{\Delta} \Delta t_{m-1}\\
    + \nabla_xb\left(t_{m-1}, X^{\Delta}_{m-1}\right) D_n X_{m-1}^{\Delta} \Delta W_{m-1},  \quad  0 \leq n < m \leq N.
    \end{cases}
     \label{eq9}
\end{equation}
and
\begin{equation}
    \begin{split}
    D_n Y_n^{\Delta} & = D_n Y_{n+1}^{\Delta} + f_D\left(t_n, \mathbf{X}_n^{\Delta}, \mathbf{D}_n \mathbf{X}_n^{\Delta}\right)\,\Delta t_n - D_n Z_n^{\Delta} \Delta W_n,
    \end{split}
    \label{eq10}
\end{equation}
where $\mathbf{D}_n\mathbf{X}_n^{\Delta} := \left( D_n X_n^{\Delta}, D_n Y_n^{\Delta}, D_n Z_n^{\Delta}\right)$. Note that $D_n Y_n^{\Delta} = Z_n^{\Delta}$ due to~\eqref{eq8}. Moreover, using the Malliavin chain rule~\cite{nualart2006malliavin} and the Feynman-Kac relation~\eqref{eq3}, we have that
$$\quad D_n Y_{n+1}^{\Delta} = \nabla_x u\left(t_{n+1}, X_{n+1}^{\Delta}\right)D_n X_{n+1}^{\Delta}, \quad D_n Z_n^{\Delta} =: \gamma\left(t_{n}, X_{n}^{\Delta}\right)D_n X_{n}^{\Delta},$$
where $\gamma: [0, T] \times \mathbb{R}^{d} \to \mathbb{R}^{d \times d} $ is the Jacobian matrix of $\nabla_x u(t, x) b(t,x)$. 

After discretizing the integrals, our scheme is made fully implementable by approximating the unknown processes $\left(Y_n^{\Delta}, Z_n^{\Delta}, \Gamma_n^{\Delta}\right)$ in the discrete BSDE system~\eqref{eq5} and~\eqref{eq10} using three DNNs for $n=0,1,\ldots,N$, where $\Gamma_n^{\Delta} := \gamma(t_n, X_n^{\Delta})$. We refer to our scheme as differential LDBSDE (DLDBSDE), which works as follows:
\begin{itemize}
    \item Generate approximations $X^{\Delta}_{n+1}$ for $n = 0, 1, \ldots, N-1$ using~\eqref{eq4} and its discrete Malliavin derivative $D_nX_n^{\Delta}, D_nX_{n+1}^{\Delta}$ using~\eqref{eq9}.
    \item At each discrete time point $t_n$, $n = 0, 1, \ldots, N$, 
    use DNNs $\phi^y(\cdot; \theta^y): \mathbb{R}^{1+d} \to \mathbb{R}$, $\phi^z(\cdot; \theta^z): \mathbb{R}^{1+d} \to \mathbb{R}^{1 \times d}$ and $\phi^{\gamma}(\cdot; \theta^{\gamma}): \mathbb{R}^{1+d} \to \mathbb{R}^{d \times d}$ to approximate the
    discrete processes $\left(Y^{\Delta}_n, Z^{\Delta}_n, \Gamma^{\Delta}_n\right)$, respectively, where the input vector of the network is the time value $t_n \in \mathbb{R}_{+}$ and the Markovian process $X_n^{\Delta} \in \mathbb{R}^d$, namely
    $$Y_n^{\Delta,\theta} = \phi^y(t_n, X_n^{\Delta};\theta^y), \quad Z_n^{\Delta,\theta} = \phi^z(t_n, X_n^{\Delta};\theta^z), \quad \Gamma_n^{\Delta,\theta} = \phi^{\gamma}(t_n, X_n^{\Delta};\theta^{\gamma}).$$
    \item Train the parameters $\theta=\left( \theta^y, \theta^z, \theta^{\gamma} \right)$ using a global differential loss type function including local losses such that the dynamics of discretized BSDE system~\eqref{eq5} and~\eqref{eq10} are satisfied at each time step, namely
    \begin{equation}
        \begin{aligned}            
            \mathbf{L}^{\Delta}\left( \theta \right)& := \omega_1 \mathbf{L}^{y,\Delta}\left( \theta \right) + \omega_2 \mathbf{L}^{z,\Delta}\left( \theta \right)\\
            \mathbf{L}^{y,\Delta}\left( \theta \right)& = \mathbb{E}\left[\sum_{n=0}^{N-1} \left\vert Y^{\Delta, \theta}_{n+1} - Y_n^{\Delta, \theta} + f\left(t_n, \mathbf{X}^{\Delta, \theta}_n\right) \Delta t -  Z_n^{\Delta,\theta} \Delta W_n  \right\vert^2 + \left\vert Y_N^{\Delta, \theta} - g\left(X_N^{\Delta}\right)\right\vert^2\right]\\
            \mathbf{L}^{z,\Delta}\left( \theta_n \right) & := \mathbb{E}\left[ \sum_{n=0}^{N-1}\left\vert D_n Y^{\Delta, \theta}_{n+1} - Z_n^{\Delta,\theta} + f_D\left(t_n, \mathbf{X}^{\Delta, \theta}_n, \mathbf{D}_n\mathbf{X}_n^{\Delta, \theta}\right)\Delta t - \Gamma_n^{\Delta,\theta} D_n X_n^{\Delta} \Delta W_n
        \right\vert^2 \right.\\
         & \quad\quad\quad \left. + \left\vert Z_N^{\Delta, \theta} - g_x\left(X_N^{\Delta}\right)b\left( t_{N}, X_{N}^{\Delta} \right)\right\vert^2 \vphantom{\sum_{n=0}^{N-1}} \right],
        \end{aligned}
        \label{eq11}
    \end{equation}
     where $D_n Y^{\Delta, \theta}_{n+1} = Z^{\Delta, \theta}_{n+1} b^{-1}\left( t_{n+1}, X_{n+1}^{\Delta} \right) D_n X^{\Delta}_{n+1}$, $\omega_1, \omega_2 \in [0, 1]$, $\omega_1 +\omega_2 = 1$,  and for notational convinience $\mathbf{X}^{\Delta, \theta}_n:=\left( X_n^{\Delta}, Y_n^{\Delta, \theta}, Z_n^{\Delta, \theta} \right)$ and $\mathbf{D}_n\mathbf{X}^{\Delta, \theta}_n:=\left( D_n X_n^{\Delta}, Z_n^{\Delta,\theta}, \Gamma_n^{\Delta,\theta} D_n X_n^{\Delta} \right)$. 
     \item Approximate the optimal parameters $\theta^* \in \argmin_{\theta \in \Theta} \mathbf{L}^{\Delta}\left( \theta \right)$ using a SGD method and receive the final estimated parameters $\hat{\theta}$. Set the final approximation of $\left( Y_n^{\Delta}, Z_n^{\Delta}, \Gamma_n^{\Delta} \right)$ as $\left( Y_n^{\Delta, \hat{\theta}}, Z_n^{\Delta, \hat{\theta}}, \Gamma_n^{\Delta, \hat{\theta}}\right)$ for $n=0,1,\ldots,N$. 
\end{itemize}
Note that LDBSDE scheme can be considered as a special case of our scheme by choosing $\omega_1 = 1$ and $\omega_2 = 0$, and using one DNN for $Y$ and AD for the processes $Z$ and $\Gamma$. For our scheme, $\omega_1 = \frac{1}{d+1}$ and $\omega_2 = \frac{d}{d+1}$ are considered due to corresponding dimensionality of the processes $Y$ and $Z$, a practice used in differential deep learning~\cite{huge2020differential}. The difference between the DLDBSDE scheme and the backward differential deep learning scheme in~\cite{kapllani2024backward} is outlined in Remark~\ref{r1}, and its convergence analysis in Remark~\ref{r2}.
\begin{remark}
The DLDBSDE scheme significantly differs from the backward differential deep learning scheme presented in~\cite{kapllani2024backward}. While our previous work focused on developing a backward type method to provide highly accurate first- and second-order gradient approximations, DLDBSDE scheme introduces a novel forward one. Unlike the backward approach, which relied on local optimization, our current method is based on global optimization. Additionally, the computation time of the scheme in~\cite{kapllani2024backward} is comparable to its deep learning counterpart~\cite{hure2020deep} when including the computation of 
$\Gamma$ at each optimization step. In contrast, the computation time of the DLDBSDE scheme is significantly shorter compared to the its deep learning counterpart, namely the LDBSDE scheme, as we demonstrate in the numerical section.
\label{r1}
\end{remark}
\begin{remark}
The convergence of the LDBSDE scheme (and a Milstein-version of the scheme) is discussed in~\cite{knochenhauer2021convergence}, providing an a posteriori error estimation similar to~\cite{han2020convergence} for the DBSDE scheme. The authors demonstrate that the error of the LDBSDE scheme is bounded by its respective loss function~\eqref{eq6}, and the loss functional converges sufficiently fast to zero, ensuring that the error of the scheme vanishes in the limit. This result is achievable through the universal approximation theorem~\cite{hornik1989multilayer,cybenko1989approximation} of neural networks. An a posteriori error analysis for the DLDBSDE scheme can be conducted by following the methodology in~\cite{knochenhauer2021convergence} and our work in~\cite{kapllani2024backward}. The latter includes  the additional assumptions needed for ensuring the boundedness of the Malliavin derivatives. It also provides the extra steps required to address the discretization error introduced by the Euler-Maruyama method in~\eqref{eq9}-\eqref{eq10} and the model/approximation error from the DNNs $\left(\phi^z, \phi^{\gamma}\right)$ associated with the second term in the loss function~\eqref{eq11}. This is part of our ongoing research.
\label{r2}
\end{remark}

\section{Numerical results}
\label{sec4}
In this section, we demonstrate the improved performance of the DLDBSDE scheme compared to the LDBSDE scheme. We demonstrate this improvement in approximating the solution, its gradient and the Hessian matrix. Given the importance of accurate gradient approximations in finance, particularly in delta- and $\Gamma$-hedging strategies for option contracts, we also focus on option pricing examples. All the experiments below were run in PYTHON using TensorFlow on the PLEIADES cluster (no parallelization), which consists of 268 workernodes and additionally 5 GPU nodes with 8 NVidia HGX A100 GPUs (128 cores each, 2 TB memory, 16 GB per thread). We run the algorithms on the GPU nodes. For more information, see~PLEIADES documentation\footnote{\url{https://pleiadesbuw.github.io/PleiadesUserDocumentation/}}.

In each example, we use the same hyperparameters for both our scheme and the LDBSDE scheme to ensure a fair comparison. For the DNNs, we choose $L = 4$ hidden layers and $\eta = 100 + d$ neurons per hidden layer. The input is normalized based on the true moments, excluding the normalization at discrete time point $t_0$ due to zero standard deviation. We employ a hyperbolic tangent activation $\tanh(\cdot)$ on each hidden layer. The Adam optimizer is used as an SGD algorithm with a stepwise learning rate decay approach. The total number of optimization steps is $\Kf = 60000$, where a batch size of $B = 128$ is considered for each step $\kappa$, and the learning rate $\alpha$ is adjusted as follows 
\begin{equation*}
  \alpha_{\kappa}=\begin{cases}
    \num{1e-3}, & \text{for $1 \leq \kappa \leq 20000$},\\
    \num{3e-4}, & \text{for $20000 < \kappa \leq 30000$},\\
    \num{1e-4}, & \text{for $30000 < \kappa \leq 4000$},\\
    \num{3e-5}, & \text{for $40000 < \kappa \leq 50000$},\\
    \num{1e-5}, & \text{for $50000 < \kappa \leq \Kf$}.\\
    \end{cases}
\end{equation*}
We define the mean squared errors (MSEs) as performance metrics for a sample of size $B$:
\begin{equation*}
    \tilde{\varepsilon}^{y}_n := \frac{1}{B}\sum_{j=1}^B \left\vert Y_{n, j} - Y_{n, j}^{\Delta, \hat{\theta}} \right\vert^2, \quad \tilde{\varepsilon}^{z}_n := \frac{1}{B}\sum_{j=1}^B \left\vert Z_{n, j} - Z_{n, j}^{\Delta, \hat{\theta}} \right\vert^2, \quad \tilde{\varepsilon}^{\gamma}_n := \frac{1}{B}\sum_{j=1}^B \left\vert \Gamma_{n, j} - \Gamma_{n, j}^{\Delta, \hat{\theta}} \right\vert^2,
\end{equation*}
for each process. To account the stochasticity of the underlying Brownian motion and the Adam optimizer, we conduct $Q = 10$ independent runs (trainings) of the algorithms. We then define, e.g.,
\begin{equation*}
    \overline{{\tilde{\varepsilon}}}^{y}_n := \frac{1}{Q} \sum_{q=1}^Q \tilde{\epsilon}^{y}_{n,q},
\end{equation*}
as the mean MSE for the process $Y$, and similarly for the other processes. As a relative measure of the MSE, we consider,
\begin{equation*}
    \tilde{\varepsilon}^{y, r}_n := \frac{1}{B}\sum_{j=1}^B \frac{\left\vert Y_{n, j} - Y_{n, j}^{\Delta, \hat{\theta}} \right\vert^2}{\left\vert Y_{n, j} \right\vert^2},
\end{equation*} 
for the process $Y$, and similarly for the other processes. We select a testing sample of size $B = 1024.$ The computation time (runtime) in seconds for one run of the algorithms is denoted as $\tau$, and the average computation time over $Q=10$ runs as $\overline{{\tau}}: = \frac{1}{Q} \sum_{q=1}^Q \tau_q.$

\subsection{The simple bounded BSDE}
We start with the simple bounded BSDE studied in~\cite{hure2020deep,teng2022gradient,chassagneux2023learning,kapllani2024deep}.
\label{subsec42}
\begin{example}
The high-dimensional BSDE given in~\cite{hure2020deep} reads
\begin{equation*}
    \begin{split}
        \left\{
            \begin{array}{rcl}
                dX_t & = & a \,dt + b \, dW_t,\\
                X_0 & = & x_0,\\
                -dY_t & = & \left(\left( \cos\left(\sum_{k=1}^{d}X_t^k\right)  + 0.2 \sin\left(\sum_{k=1}^{d}X_t^k\right)\right)\exp\left(\frac{T-t}{2}\right) \right. \\
                & & \left. -\frac{1}{2} \left( \sin\left(\sum_{k=1}^{d}X_t^k\right) \cos\left(\sum_{k=1}^{d}X_t^k\right) \exp\left( T-t\right) \right)^2  + \frac{1}{2d}\left(Y_t \sum_{k=1}^{d}Z_t^k\right)^2 \vphantom{\cos} \right)\,dt\\
                & & -Z_t \,dW_t,\\  
                Y_T & = & \cos\left(\sum_{k=1}^{d}X_t^k\right).
            \end{array}
        \right. 
    \end{split}
\end{equation*}
\label{ex1}
\end{example}
The analytical solution is given by
\begin{equation*}
    \begin{split}
        \left\{
            \begin{array}{rcl}
                Y_t & = & \exp\left(\frac{T-t}{2}\right) \cos\left( \sum_{k=1}^{d}X_t^k \right),\\
   	   	    	Z_t & = &-b \exp\left(\frac{T-t}{2}\right) \sin\left(\sum_{k=1}^{d}X_t^k \right)\mathbf{1}_{1,d}.
            \end{array}
        \right. 
    \end{split}
\end{equation*}
Note that the analytical solution $\Gamma_t = \nabla_x \left( \nabla_x u\left(t, X_t\right) b\left(t, X_t\right) \right)$ is calculated by using AD (similarly for the following examples). We choose $a = \frac{0.2}{d}$, $b = \frac{1}{\sqrt{d}}$ and $x_0 = \mathbf{1}_d$. This example is very interesting as for $d=1$ and $T=2$, the DBSDE scheme~\cite{weinan2017deep} diverges, while the LDBSDE scheme converges to an approximation far from the exact solution, namely to a poor local minima, see~\cite{kapllani2024deep}. Hence, we test if for each process, the DLDBSDE scheme can converge to a better local minima compared to the LDBSDE scheme. In Table~\ref{tab1}, we report the mean relative MSE values at $t_0$
for $\left(Y_0, Z_0, \Gamma_0\right)$ from the LDBSDE and DLDBSDE schemes, their average runtime (in seconds) and the empirical convergence rates $\beta$ using $N \in \{4, 16, 64, 256\}$. The STD of the relative MSE values at $t_0$ is given in the brackets. Note that for $N=256$ the approximations from the LDBSDE are not available, because the scheduled scripts in the GPU nodes of PLEIADES cluster have a time limit of 3 days.
\begin{table}[htb!]
{\footnotesize
\begin{center}
  \resizebox{\textwidth}{!}{\begin{tabular}{| c | c | c | c | c | c |}
  \hline
    \multirow{3}{*}{Metric} & N = 4 & N = 16 & N = 64 & N = 256 & \multirow{3}{*}{$\beta$}\\ \
    & LDBSDE & LDBSDE & LDBSDE & LDBSDE & \\ 
    & DLDBSDE & DLDBSDE & DLDBSDE & DLDBSDE & \\ \hline
    \multirow{2}{*}{$\overline{{\tilde{\varepsilon}}}^{y, r}_0$} & $\num{3.28e+00}$ $(\num{4.14e-02})$ & $\num{4.49e-01}$ $(\num{5.57e-03})$ & $\num{1.88e-01}$ $(\num{1.83e-03})$ & NA & $-$ \\
    & $\num{6.35e+00}$ $(\num{5.87e-03})$ & $\num{2.35e-01}$ $(\num{1.33e-03})$ & $\num{1.04e-02}$ $(\num{4.69e-04})$ & $\num{7.09e-03}$ $(\num{2.56e-04})$ & $1.70$ \\ \hline    
    \multirow{2}{*}{$\overline{{\tilde{\varepsilon}}}^{z, r}_0$} & $\num{2.55e-01}$ $(\num{7.35e-03})$ & $\num{5.83e-02}$ $(\num{4.19e-04})$ & $\num{3.38e-02}$ $(\num{5.33e-04})$ & NA & $-$ \\
    &  $\num{1.75e+00}$ $(\num{1.66e-03})$ & $\num{2.10e-01}$ $(\num{9.37e-04})$ & $\num{2.35e-03}$ $(\num{7.68e-05})$ & $\num{7.62e-03}$ $(\num{6.41e-04})$ & $1.50$ \\ \hline    
    \multirow{2}{*}{$\overline{{\tilde{\varepsilon}}}^{\gamma, r}_0$} & $\num{2.63e+00}$ $(\num{3.32e-01})$ & $\num{7.02e-01}$ $(\num{1.32e-01})$ & $\num{3.96e+00}$ $(\num{2.55e+00})$ & NA & $-$ \\
    & $\num{5.60e-02}$ $(\num{8.13e-04})$ & $\num{9.76e-03}$ $(\num{7.67e-04})$ & $\num{6.84e-02}$ $(\num{1.40e-03})$ & $\num{3.76e-01}$ $(\num{2.29e-02})$ & $-0.55$ \\ \hline    
    \multirow{2}{*}{$\overline{\tau}$} & $\num{7.15e+03}$ & $\num{2.33e+04}$ & $\num{9.08e+04}$ & NA & \\  
    & $\num{6.38e+02}$ & $\num{1.56e+03}$ & $\num{5.14e+03}$ & $\num{2.03e+04}$ &  \\ \hline    
   \end{tabular}}
  \end{center}
\caption{Mean relative MSE values, empirical convergence rates of $\left(Y_0, Z_0, \Gamma_0 \right)$ from LDBSDE and DLDBSDE schemes and their average runtimes in Example~\ref{ex1} for $d=1$, $T=2$ and $N \in \{4, 16, 64, 256\}$. The STD of the relative MSE values at $t_0$ is given in the brackets. The approximations for $N=256$ from the LDBSDE are not available (NA) due to large computation time (more than 3 days).}
\label{tab1} 
}
\end{table}
Our scheme provides the smallest mean relative MSE for $Y_0$. This is also observed for $Z_0$, and especially for $\Gamma_0$. Additionally, the computational cost of the DLDBSDE algorithm is significantly lower than that of the LDBSDE scheme. The empirical convergence rates in this and subsequent examples can be improved by reducing the optimization and model errors. This can be achieved by increasing the optimization steps $\Kf$, the batch size $B$, and the number of hidden neurons $\eta$ or hidden layers $L$.

\begin{figure}[tbhp!]
	\centering
	\subfloat[Process $Y$.]{
        \pgfplotstableread{"Figures/Example1/d1/vareps_Y.dat"}{\table}        
        \begin{tikzpicture} 
            \begin{axis}[
                xmin = 0, xmax = 2.0,
                ymin = 5e-5, ymax = 5e-1,
                xtick distance = 0.4,
                ymode=log, 
                grid = both,
                width = 0.33\textwidth,
                height = 0.35\textwidth,
                xlabel = {$t_n$},
                legend cell align = {left},
                legend pos = south west,
                legend style={nodes={scale=0.5, transform shape}}]
                \addplot[smooth, ultra thick, blue, dashed] table [x = 0, y =1] {\table}; 
                \addplot[smooth, ultra thick, green, solid] table [x = 0, y =4] {\table}; 
                \addplot [name path=upper,draw=none] table[x=0,y=2] {\table};
                \addplot [name path=lower,draw=none] table[x=0,y=3] {\table};
                \addplot [fill=blue!40] fill between[of=upper and lower];           
                \addplot [name path=upper,draw=none] table[x=0,y=5] {\table};
                \addplot [name path=lower,draw=none] table[x=0,y=6] {\table};
                \addplot [fill=green!40] fill between[of=upper and lower];                
                \legend{
                   $\overline{{\tilde{\varepsilon}}}^{y}_n$-LDBSDE,
                   $\overline{{\tilde{\varepsilon}}}^{y}_n$-DLDBSDE
                }
            \end{axis}
            \end{tikzpicture}
		\label{fig1a}
	} 
	\subfloat[Process $Z$.]{ \pgfplotstableread{"Figures/Example1/d1/vareps_Z.dat"}{\table}
        \begin{tikzpicture} 
            \begin{axis}[
                xmin = 0, xmax = 2.0,
                ymin = 1e-4, ymax = 2e-1,
                xtick distance = 0.4,
                ymode=log,
                grid = both,
                width = 0.33\textwidth,
                height = 0.35\textwidth,
                xlabel = {$t_n$},
                legend cell align = {left},
                legend pos = south west,
                legend style={nodes={scale=0.5, transform shape}}]
                \addplot[smooth, ultra thick, blue, dashed] table [x = 0, y =1] {\table}; 
                \addplot[smooth, ultra thick, green, solid] table [x = 0, y =4] {\table}; 
                \addplot [name path=upper,draw=none] table[x=0,y=2] {\table};
                \addplot [name path=lower,draw=none] table[x=0,y=3] {\table};
                \addplot [fill=blue!40] fill between[of=upper and lower];           
                \addplot [name path=upper,draw=none] table[x=0,y=5] {\table};
                \addplot [name path=lower,draw=none] table[x=0,y=6] {\table};
                \addplot [fill=green!40] fill between[of=upper and lower];                
                \legend{
                   $\overline{{\tilde{\varepsilon}}}^{z}_n$-LDBSDE,
                   $\overline{{\tilde{\varepsilon}}}^{z}_n$-DLDBSDE
                }
            \end{axis}
            \end{tikzpicture}
		\label{fig1b}
	} 
        \subfloat[Process $\Gamma$.]{ \pgfplotstableread{"Figures/Example1/d1/vareps_Gamma.dat"}{\table}
        \begin{tikzpicture} 
            \begin{axis}[
                xmin = 0, xmax = 2.0,
                ymin = 8e-4, ymax = 1e+1,
                xtick distance = 0.4,
                ymode=log,
                grid = both,
                width = 0.33\textwidth,
                height = 0.35\textwidth,
                xlabel = {$t_n$},
                legend cell align = {left},
                legend pos = south west,
                legend style={nodes={scale=0.5, transform shape}}]
                \addplot[smooth, ultra thick, blue, dashed] table [x = 0, y =1] {\table}; 
                \addplot[smooth, ultra thick, green, solid] table [x = 0, y =4] {\table}; 
                \addplot [name path=upper,draw=none] table[x=0,y=2] {\table};
                \addplot [name path=lower,draw=none] table[x=0,y=3] {\table};
                \addplot [fill=blue!40] fill between[of=upper and lower];           
                \addplot [name path=upper,draw=none] table[x=0,y=5] {\table};
                \addplot [name path=lower,draw=none] table[x=0,y=6] {\table};
                \addplot [fill=green!40] fill between[of=upper and lower];                
                \legend{
                   $\overline{{\tilde{\varepsilon}}}^{\gamma}_n$-LDBSDE,
                   $\overline{{\tilde{\varepsilon}}}^{\gamma}_n$-DLDBSDE
                }
            \end{axis}
            \end{tikzpicture}
		\label{fig1c}
	} 
    \caption{Mean MSE values of the processes $\left(Y, Z, \Gamma \right)$ from LDBSDE and DLDBSDE schemes over the discrete time points $\{t_n\}_{n=0}^{N}$ using the testing sample in Example~\ref{ex1}, for $d=1$, $T=2$ and $N = 64$. The STD of MSE values is given in the shaded area.}
\label{fig1}
\end{figure}
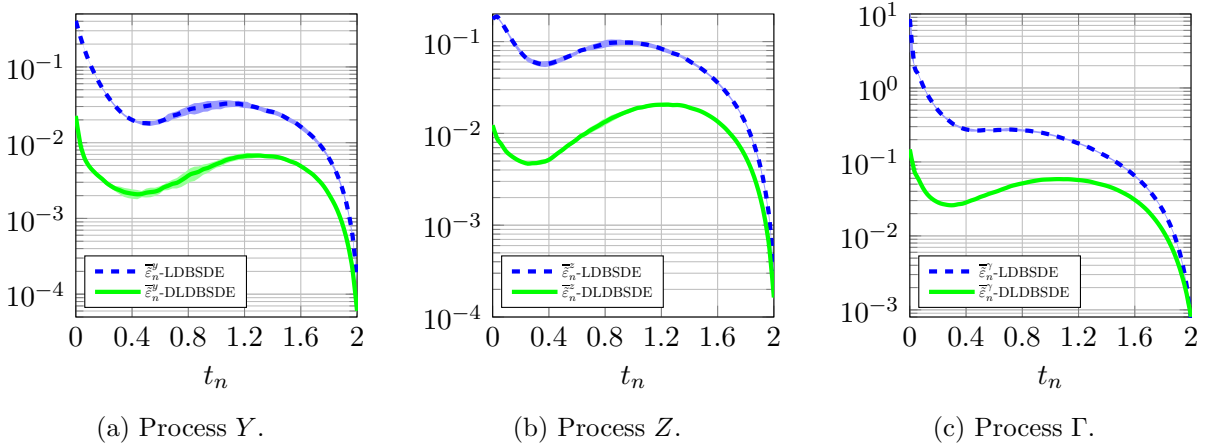

To provide a comparison of the approximation of each process using the entire testing sample over discrete domain $\Delta$, we display in Figure~\ref{fig1} the mean MSE values for $\left( Y_n^{\Delta, \hat{\theta}}, Z_n^{\Delta, \hat{\theta}}, \Gamma_n^{\Delta, \hat{\theta}} \right)$, $n = 0, 1, \ldots, N$, using the testing sample of size $B=1024$ and $N = 64$. The STD of the MSE values is visualized in the shaded area.
Our scheme outperforms the LDBSDE scheme in approximating each process across the discrete time points $\{t_n\}_{n=0}^{N}$. Figure~\ref{fig1c} shows a substantial improvement in approximating the process $\Gamma$.

Next, we increase the dimension to $d=50$ and choose $T=0.5$. In Table~\ref{tab2}, we report the mean relative MSE values at $t_0$ for each process, the algorithm average runtime and the empirical convergence rates using $N \in \{ 2, 8, 32, 64\}$. The STD of the relative MSE values at $t_0$ is given in the brackets.
\begin{table}[htb!]
{\footnotesize
\begin{center}
  \resizebox{\textwidth}{!}{\begin{tabular}{| c | c | c | c | c | c |}
  \hline
    \multirow{3}{*}{Metric} & N = 2 & N = 8 & N = 32 & N = 64 & \multirow{3}{*}{$\beta$} \\ \
    & LDBSDE & LDBSDE & LDBSDE & LDBSDE & \\ 
    & DLDBSDE & DLDBSDE & DLDBSDE & DLDBSDE & \\ \hline
    \multirow{2}{*}{$\overline{{\tilde{\varepsilon}}}^{y, r}_0$} & $\num{4.14e-04}$ $(\num{1.98e-04})$ & $\num{3.69e-05}$ $(\num{5.36e-05})$ & $\num{2.24e-04}$ $(\num{5.56e-05})$ & $\num{2.69e-04}$ $(\num{7.74e-05})$ & $-0.01$\\
    & $\num{5.39e-03}$ $(\num{2.23e-04})$ & $\num{4.72e-04}$ $(\num{4.30e-05})$ & $\num{1.19e-04}$ $(\num{2.36e-05})$ & $\num{1.08e-04}$ $(\num{2.24e-05})$ & $1.15$ \\ \hline    
    \multirow{2}{*}{$\overline{{\tilde{\varepsilon}}}^{z, r}_0$} & $\num{2.24e-01}$ $(\num{3.89e-03})$ & $\num{5.24e-02}$ $(\num{4.55e-03})$ & $\num{9.17e-01}$ $(\num{3.93e-02})$ & $\num{9.80e-01}$ $(\num{1.01e-02})$ & $-0.60$ \\
    & $\num{6.53e-03}$ $(\num{1.14e-03})$ & $\num{5.36e-04}$ $(\num{2.54e-04})$ & $\num{1.34e-04}$ $(\num{8.38e-05})$ & $\num{1.29e-04}$ $(\num{1.15e-04})$ & $1.16$ \\ \hline    
    \multirow{2}{*}{$\overline{{\tilde{\varepsilon}}}^{\gamma, r}_0$} & $\num{9.86e-01}$ $(\num{1.05e-02})$ & $\num{1.04e+00}$ $(\num{8.73e-03})$ & $\num{1.00e+00}$ $(\num{7.24e-04})$ & $\num{1.00e+00}$ $(\num{1.97e-04})$ & $-0.01$ \\
    & $\num{4.56e-02}$ $(\num{1.53e-03})$ & $\num{1.91e-02}$ $(\num{4.71e-04})$ & $\num{1.73e-02}$ $(\num{5.63e-04})$ & $\num{2.85e-02}$ $(\num{3.03e-03})$ & $0.16$ \\ \hline    
    \multirow{2}{*}{$\overline{\tau}$} & $\num{3.66e+03}$ & $\num{1.09e+04}$ & $\num{4.07e+04}$ & $\num{8.11e+04}$ & \\  
    & $\num{1.79e+03}$ & $\num{5.72e+03}$ & $\num{2.15e+04}$ & $\num{4.26e+04}$ & \\ \hline    
   \end{tabular}}
  \end{center}
\caption{Mean relative MSE values, empirical convergence rates of $\left(Y_0, Z_0, \Gamma_0 \right)$ from LDBSDE and DLDBSDE schemes and their average runtimes in Example~\ref{ex1} for $d=50$, $T=0.5$ and $N \in \{2, 8, 32, 64\}$. The STD of the relative MSE values at $t_0$ is given in the brackets.}
\label{tab2} 
}
\end{table}
The approximations for $Y_0$ are comparable in both the schemes. However, for $\left(Z_0, \Gamma_0\right)$, the approximations are significantly more accurate with the DLDBSDE scheme compared to the LDBSDE scheme, and this increased accuracy is achieved with considerably less computation time. The mean MSE values over the entire discrete domain $\Delta$ for the testing sample are given in Figure~\ref{fig2}, where the STD of the MSE values is given in the shaded area.
\begin{figure}[tbhp!]
	\centering
	\subfloat[Process $Y$.]{
        \pgfplotstableread{"Figures/Example1/d50/vareps_Y.dat"}{\table}        
        \begin{tikzpicture} 
            \begin{axis}[
                xmin = 0, xmax = 0.5,
                ymin = 3e-6, ymax = 7e-3,
                xtick distance = 0.1,
                ymode=log, 
                grid = both,
                width = 0.33\textwidth,
                height = 0.35\textwidth,
                xlabel = {$t_n$},
                legend cell align = {left},
                legend pos = south west,
                legend style={nodes={scale=0.5, transform shape}}]
                \addplot[smooth, ultra thick, blue, dashed] table [x = 0, y =1] {\table}; 
                \addplot[smooth, ultra thick, green, solid] table [x = 0, y =4] {\table}; 
                \addplot [name path=upper,draw=none] table[x=0,y=2] {\table};
                \addplot [name path=lower,draw=none] table[x=0,y=3] {\table};
                \addplot [fill=blue!40] fill between[of=upper and lower];           
                \addplot [name path=upper,draw=none] table[x=0,y=5] {\table};
                \addplot [name path=lower,draw=none] table[x=0,y=6] {\table};
                \addplot [fill=green!40] fill between[of=upper and lower];                
                \legend{
                   $\overline{{\tilde{\varepsilon}}}^{y}_n$-LDBSDE,
                   $\overline{{\tilde{\varepsilon}}}^{y}_n$-DLDBSDE
                }
            \end{axis}
            \end{tikzpicture} 
		\label{fig2a}
	} 
	\subfloat[Process $Z$.]{ \pgfplotstableread{"Figures/Example1/d50/vareps_Z.dat"}{\table}
        \begin{tikzpicture} 
            \begin{axis}[
                xmin = 0, xmax = 0.5,
                ymin = 2e-7, ymax = 1e-1,
                xtick distance = 0.1,
                ymode=log,
                grid = both,
                width = 0.33\textwidth,
                height = 0.35\textwidth,
                xlabel = {$t_n$},
                legend cell align = {left},
                legend pos = south west,
                legend style={nodes={scale=0.5, transform shape}}]
                \addplot[smooth, ultra thick, blue, dashed] table [x = 0, y =1] {\table}; 
                \addplot[smooth, ultra thick, green, solid] table [x = 0, y =4] {\table}; 
                \addplot [name path=upper,draw=none] table[x=0,y=2] {\table};
                \addplot [name path=lower,draw=none] table[x=0,y=3] {\table};
                \addplot [fill=blue!40] fill between[of=upper and lower];           
                \addplot [name path=upper,draw=none] table[x=0,y=5] {\table};
                \addplot [name path=lower,draw=none] table[x=0,y=6] {\table};
                \addplot [fill=green!40] fill between[of=upper and lower];                
                \legend{
                   $\overline{{\tilde{\varepsilon}}}^{z}_n$-LDBSDE,
                   $\overline{{\tilde{\varepsilon}}}^{z}_n$-DLDBSDE
                }
            \end{axis}
            \end{tikzpicture}
		\label{fig2b}
	} 
        \subfloat[Process $\Gamma$.]{ \pgfplotstableread{"Figures/Example1/d50/vareps_Gamma.dat"}{\table}
        \begin{tikzpicture} 
            \begin{axis}[
                xmin = 0, xmax = 0.5,
                ymin = 2e-3, ymax = 1e+2,
                xtick distance = 0.1,
                ymode=log,
                grid = both,
                width = 0.33\textwidth,
                height = 0.35\textwidth,
                xlabel = {$t_n$},
                legend cell align = {left},
                legend pos = south west,
                legend style={nodes={scale=0.5, transform shape}}]
                \addplot[smooth, ultra thick, blue, dashed] table [x = 0, y =1] {\table}; 
                \addplot[smooth, ultra thick, green, solid] table [x = 0, y =4] {\table}; 
                \addplot [name path=upper,draw=none] table[x=0,y=2] {\table};
                \addplot [name path=lower,draw=none] table[x=0,y=3] {\table};
                \addplot [fill=blue!40] fill between[of=upper and lower];           
                \addplot [name path=upper,draw=none] table[x=0,y=5] {\table};
                \addplot [name path=lower,draw=none] table[x=0,y=6] {\table};
                \addplot [fill=green!40] fill between[of=upper and lower];                
                \legend{
                   $\overline{{\tilde{\varepsilon}}}^{\gamma}_n$-LDBSDE,
                   $\overline{{\tilde{\varepsilon}}}^{\gamma}_n$-DLDBSDE
                }
            \end{axis}
            \end{tikzpicture}
		\label{fig2c}
	} 
    \caption{Mean MSE values of the processes $\left(Y, Z, \Gamma \right)$ from LDBSDE and DLDBSDE schemes over the discrete time points $\{t_n\}_{n=0}^{N}$ using the testing sample in Example~\ref{ex1}, for $d=50$, $T=0.5$ and $N = 64$. The STD of MSE values is given in the shaded area.}
\label{fig2}
\end{figure}
As for $d=1$, our method shows better approximations of each process on the entire time domain compared to the LDBSDE scheme.

\subsection{The Black-Scholes BSDE}
\label{subsec43}
We now consider a linear option pricing example, the Black-Scholes BSDE, which is used for pricing of European options. 
\begin{example}
The high-dimenisonal Black-Scholes BSDE is given as follows~\cite{zhang2013sparse}
\begin{equation*}
    \begin{split}
    \left\{
        \begin{array}{rcl}
             dX_t^k &=&  \left(a_k - \delta_k\right) X_t^k\,dt + b_k X_t^k \,dW_t^k, \\
             X_0^k &=& x_0^k, \quad k = 1,\ldots, d,\\  
           -dY_t &=& - \left( R Y_t + \sum_{k=1}^d  \frac{ a_k - R + \delta_k}{b_k} Z_t^k\right)\,dt- Z_t \,dW_t,\\ 
   		   	   Y_T &=& \left(\prod_{k=1}^d \left(X_T^k\right)^{c_k}-K\right)^+,
        \end{array}
    \right. \\ 
    \end{split}
\end{equation*}
\label{ex2}
\end{example}
where $c_k >0$ and $\sum_{k=1}^d c_k = 1$. Note that $a_k$ represents the return rate of the stock $X_t^k$, $b_k$ the volatility of the stock returns, $\delta_k$ is its dividend rate, and $x_0^k$ is the price of the stock at $t =0$. Moreover, $T$ denotes the maturity of the option contract, while $K$ represents the contract's strike price. Finally, $R$ corresponds to the risk-free interest rate. The analytic solution (the option price $Y_t$ and delta-hedging strategy $Z_t$) is given by
\begin{equation}
 \begin{split}
    \left\{
        \begin{array}{rcl}
            Y_t &=& u(t, X_t) = \exp\left(-\check{\delta} \left(T-t\right)\right) \prod_{k=1}^d \left(X_t^k\right)^{c_k} \Phi \left(\check{d}_1\right)-\exp\left(-R\left(T-t\right)\right) K \Phi \left(\check{d}_2\right),\\
           Z_t^k &=& \frac{\partial u}{\partial x_k} b_k X_t^k = c_k \exp\left(-\check{\delta} \left(T-t\right)\right) \prod_{k=1}^d \left(X_t^k\right)^{c_k} \Phi \left(\check{d}_1\right)b_k, \quad k = 1, \ldots, d,\\
           \check{d}_1 &=& \frac{\ln\left(\frac{\prod_{k=1}^d \left(X_t^k\right)^{c_k}}{K}\right) + \left( R-\check{\delta} + \frac{\check{b}^2}{2} \right) \left(T-t\right)}{\check{b} \sqrt{T-t}},\\
            \check{d}_2 &=& \check{d}_1 - \check{b}\sqrt{T-t},\\
            \check{b}^2 &=& \sum_{k=1}^d (b_k c_k)^2,\quad \check{\delta} = \sum_{k=1}^d c_k\left(\delta_k + \frac{b_k^2}{2} \right) - \frac{\check{b}^2}{2},               
        \end{array}
    \right. \\
    \end{split}
\label{eq12}
\end{equation}
where $\Phi \left(\cdot\right)$ is the standard normal cumulative distribution function. Instead of solving the above BSDE directly, we consider the transformed BSDE in the ln-domain, as it simplifies the Malliavin derivative of the forward process $X$, see~\cite{kapllani2024backward}. We choose $d=50$, $T = 0.5$, $x_0^k = 100$, $a_k = 0.05$, $b_k = 0.2$, $R = 0.03$, $c_k = \frac{1}{d}$, $\delta_k = 0$ for $k=1,\ldots,d,$ and $K=100$. The mean relative MSE values of $\left(Y_0, Z_0, \Gamma_0\right)$, the algorithm average runtime and the empirical convergence rates are reported in Table~\ref{tab3} using $N \in \{2, 8, 32, 64\}$, with the STD of the relative MSE values provided in the brackets.
\begin{table}[htb!]
{\footnotesize
\begin{center}
  \resizebox{\textwidth}{!}{\begin{tabular}{| c | c | c | c | c | c |}
  \hline
    \multirow{3}{*}{Metric} & N = 2 & N = 8 & N = 32 & N = 64 & \multirow{3}{*}{$\beta$}\\ \
    & LDBSDE & LDBSDE & LDBSDE & LDBSDE & \\ 
    & DLDBSDE & DLDBSDE & DLDBSDE & DLDBSDE & \\ \hline
    \multirow{2}{*}{$\overline{{\tilde{\varepsilon}}}^{y, r}_0$} & $\num{7.98e-03}$ $(\num{5.20e-03})$ & $\num{1.50e-02}$ $(\num{1.19e-02})$ & $\num{1.91e-02}$ $(\num{1.52e-03})$ & $\num{1.31e-02}$ $(\num{3.45e-04})$ & $-0.17$\\
    & $\num{4.31e-04}$ $(\num{1.08e-04})$ & $\num{2.06e-06}$ $(\num{1.56e-06})$ & $\num{5.45e-05}$ $(\num{1.73e-05})$ & $\num{8.81e-05}$ $(\num{2.51e-05})$ & $0.21$\\ \hline    
    \multirow{2}{*}{$\overline{{\tilde{\varepsilon}}}^{z, r}_0$} & $\num{6.57e-02}$ $(\num{1.22e-02})$ & $\num{2.51e-01}$ $(\num{1.45e-01})$ & $\num{9.87e-01}$ $(\num{1.20e-02})$ & $\num{9.96e-01}$ $(\num{5.72e-03})$ & $-0.83$\\
    & $\num{3.76e-03}$ $(\num{1.32e-04})$ & $\num{1.76e-03}$ $(\num{2.31e-04})$ & $\num{2.08e-03}$ $(\num{4.11e-04})$ & $\num{4.09e-03}$ $(\num{7.60e-04})$ & $-0.01$\\ \hline    
    \multirow{2}{*}{$\overline{{\tilde{\varepsilon}}}^{\gamma, r}_0$} & $\num{1.06e+00}$ $(\num{9.81e-02})$ & $\num{1.08e+00}$ $(\num{3.15e-02})$ & $\num{1.00e+00}$ $(\num{4.51e-04})$ & $\num{1.00e+00}$ $(\num{1.28e-04})$ & $0.02$\\
    & $\num{1.75e-02}$ $(\num{8.16e-04})$ & $\num{3.24e-02}$ $(\num{1.24e-03})$ & $\num{5.30e-02}$ $(\num{2.84e-03})$ & $\num{9.16e-02}$ $(\num{6.88e-03})$ & $-0.45$\\ \hline    
    \multirow{2}{*}{$\overline{\tau}$} & $\num{4.46e+03}$ & $\num{1.17e+04}$ & $\num{4.24e+04}$ & $\num{8.24e+04}$ & \\  
    & $\num{1.01e+03}$ & $\num{1.53e+03}$ & $\num{3.73e+03}$ & $\num{6.71e+03}$ & \\ \hline    
   \end{tabular}}
  \end{center}
\caption{Mean relative MSE values, empirical convergence rates of $\left(Y_0, Z_0, \Gamma_0 \right)$ from LDBSDE and DLDBSDE schemes and their average runtimes in Example~\ref{ex2} for $d=50$ and $N \in \{2, 8, 32, 64\}$. The STD of the relative MSE values at $t_0$ is given in the brackets.}
\label{tab3} 
}  
\end{table}
We observe that the DLDBSDE scheme significantly outperforms the LDBSDE scheme in approximating each process at $t_0$, for a shorter computation time. This improvement is evident across the entire discrete domain, as shown in Figure~\ref{fig3}, which visualizes the mean MSE for each process with the STD of the MSE values indicated by the shaded area.
\begin{figure}[tbhp!]
	\centering
	\subfloat[Process $Y$.]{
        \pgfplotstableread{"Figures/Example2/d50/vareps_Y.dat"}{\table}        
        \begin{tikzpicture} 
            \begin{axis}[
                xmin = 0, xmax = 0.5,
                ymin = 7e-5, ymax = 5e-2,
                xtick distance = 0.1,
                ymode=log, 
                grid = both,
                width = 0.33\textwidth,
                height = 0.35\textwidth,
                xlabel = {$t_n$},
                legend cell align = {left},
                legend pos = south west,
                legend style={nodes={scale=0.5, transform shape}}]
                \addplot[smooth, ultra thick, blue, dashed] table [x = 0, y =1] {\table}; 
                \addplot[smooth, ultra thick, green, solid] table [x = 0, y =4] {\table}; 
                \addplot [name path=upper,draw=none] table[x=0,y=2] {\table};
                \addplot [name path=lower,draw=none] table[x=0,y=3] {\table};
                \addplot [fill=blue!40] fill between[of=upper and lower];           
                \addplot [name path=upper,draw=none] table[x=0,y=5] {\table};
                \addplot [name path=lower,draw=none] table[x=0,y=6] {\table};
                \addplot [fill=green!40] fill between[of=upper and lower];                
                \legend{
                   $\overline{{\tilde{\varepsilon}}}^{y}_n$-LDBSDE,
                   $\overline{{\tilde{\varepsilon}}}^{y}_n$-DLDBSDE
                }
            \end{axis}
            \end{tikzpicture} 
		\label{fig3a}
	} 
	\subfloat[Process $Z$.]{ \pgfplotstableread{"Figures/Example2/d50/vareps_Z.dat"}{\table}
        \begin{tikzpicture} 
            \begin{axis}[
                xmin = 0, xmax = 0.5,
                ymin = 2e-3, ymax = 4e+0,
                xtick distance = 0.1,
                ymode=log,
                grid = both,
                width = 0.33\textwidth,
                height = 0.35\textwidth,
                xlabel = {$t_n$},
                legend cell align = {left},
                legend pos = south west,
                legend style={nodes={scale=0.5, transform shape}}]
                \addplot[smooth, ultra thick, blue, dashed] table [x = 0, y =1] {\table}; 
                \addplot[smooth, ultra thick, green, solid] table [x = 0, y =4] {\table}; 
                \addplot [name path=upper,draw=none] table[x=0,y=2] {\table};
                \addplot [name path=lower,draw=none] table[x=0,y=3] {\table};
                \addplot [fill=blue!40] fill between[of=upper and lower];           
                \addplot [name path=upper,draw=none] table[x=0,y=5] {\table};
                \addplot [name path=lower,draw=none] table[x=0,y=6] {\table};
                \addplot [fill=green!40] fill between[of=upper and lower];                
                \legend{
                   $\overline{{\tilde{\varepsilon}}}^{z}_n$-LDBSDE,
                   $\overline{{\tilde{\varepsilon}}}^{z}_n$-DLDBSDE
                }
            \end{axis}
            \end{tikzpicture}
		\label{fig3b}
	} 
        \subfloat[Process $\Gamma$.]{ \pgfplotstableread{"Figures/Example2/d50/vareps_Gamma.dat"}{\table}
        \begin{tikzpicture} 
            \begin{axis}[
                xmin = 0, xmax = 0.5,
                ymin = 1e-4, ymax = 5e-3,
                xtick distance = 0.1,
                ymode=log,
                grid = both,
                width = 0.33\textwidth,
                height = 0.35\textwidth,
                xlabel = {$t_n$},
                legend cell align = {left},
                legend pos = south west,
                legend style={nodes={scale=0.5, transform shape}}]
                \addplot[smooth, ultra thick, blue, dashed] table [x = 0, y =1] {\table}; 
                \addplot[smooth, ultra thick, green, solid] table [x = 0, y =4] {\table}; 
                \addplot [name path=upper,draw=none] table[x=0,y=2] {\table};
                \addplot [name path=lower,draw=none] table[x=0,y=3] {\table};
                \addplot [fill=blue!40] fill between[of=upper and lower];           
                \addplot [name path=upper,draw=none] table[x=0,y=5] {\table};
                \addplot [name path=lower,draw=none] table[x=0,y=6] {\table};
                \addplot [fill=green!40] fill between[of=upper and lower];                
                \legend{
                   $\overline{{\tilde{\varepsilon}}}^{\gamma}_n$-LDBSDE,
                   $\overline{{\tilde{\varepsilon}}}^{\gamma}_n$-DLDBSDE
                }
            \end{axis}
            \end{tikzpicture}
		\label{fig3c}
	} 
    \caption{Mean MSE values of the processes $\left(Y, Z, \Gamma \right)$ from LDBSDE and DLDBSDE schemes over the discrete time points $\{t_n\}_{n=0}^{N}$ using the testing sample in Example~\ref{ex2}, for $d=50$ and $N = 64$. The STD of MSE values is given in the shaded area.}
\label{fig3}
\end{figure}
Note that the accuracy of our scheme, especially for the process $\Gamma$, can be further increased by using a higher number of hidden neurons $\eta$. For $d=50$, the chosen number of hidden neurons $\eta=100+d$ is less than the output size $d\times d$ of the DNN for $\Gamma$. Therefore, increasing $\eta$ can enhance the accuracy of our scheme (provided that the optimization error is sufficiently small), particularly for $\Gamma$. This holds for the other examples as well. As already stated, the processes $\left(Z, \Gamma\right)$ play a crucial role in financial modelling, as they are related to delta- and $\Gamma$-hedging in option pricing. Hence, the improved accuracy of these processes from our scheme  compared to the LDBSDE scheme demonstrates its potential as a better tool for option pricing and hedging in high-dimensional settings.

\subsection{The Hamilton-Jacobi-Bellman equation}
\label{subsec64}
The final example is a Hamilton-Jacobi-Bellman (HJB) equation which admits a semi-explicit solution~\cite{weinan2017deep}. Solving this equation yields insights into the optimal investment strategy that maximizes the expected utility of an investor's terminal wealth. In this context, the process $Y$ represents the portfolio's wealth, while the process 
$Z$ denotes the holdings in each asset.
\begin{example}
The high-dimensional HJB BSDE reads
\begin{equation*}
    \begin{split}
        \left\{
            \begin{array}{rcl}
                dX_t & = & b \, dW_t,\\
                X_0 & = & x_0,\\
                -dY_t & = & -\sum_{k=1}^{d}\left(\frac{Z_t^k}{b}\right)^2\,dt-Z_t \,dW_t,\\  
                Y_T & = & g(X_T).
            \end{array}
        \right.
    \end{split}
\end{equation*}
\label{ex3}
\end{example}
The semi-explicit solution is given as~\cite{weinan2017deep}
$$
Y_t = u(t, X_t) = -\ln\left( \mathbb{E} \left[ \exp\left( -g( X_t +  b\left( W_T - W_t \right) ) \right) \right] \right),
$$
with $\left(Z_t, \Gamma_t\right)$ calculated using AD. As it is very time consuming to approximate highly accurate pathwise reference solutions $\left( Y_t, Z_t, \Gamma_t \right)$ for $t \in [0, T]$, we only calculate a reference solution at $t_0$. We set $T=0.5$, $d=50$, $X_0 = \mathbf{1}_d$, $b = \sqrt{0.2}$ and $g(x) = \ln\left( \frac{1}{2}\left( 1 + \left\vert x \right\vert^2 \right) \right)$ and approximate benchmark values $\left(Y_0, Z_0, \Gamma_0 \right)$ using $10^7$ Brownian motion samples and $50$ independent runs. In Table~\ref{tab4}, the mean relative MSE values of $\left(Y_0, Z_0, \Gamma_0\right)$, the algorithm average runtime and the empirical convergence rates are reported using $N \in \{2, 8, 32, 64\}$ (the STD of the relative MSE values provided in the brackets).
\begin{table}[htb!]
{\footnotesize
\begin{center}
  \resizebox{\textwidth}{!}{\begin{tabular}{| c | c | c | c | c | c |}
  \hline
    \multirow{3}{*}{Metric} & N = 2 & N = 8 & N = 32 & N = 64 & \multirow{3}{*}{$\beta$}\\ \
    & LDBSDE & LDBSDE & LDBSDE & LDBSDE & \\ 
    & DLDBSDE & DLDBSDE & DLDBSDE & DLDBSDE & \\ \hline
    \multirow{2}{*}{$\overline{{\tilde{\varepsilon}}}^{y, r}_0$} & $\num{1.29e-06}$ $(\num{1.08e-06})$ & $\num{8.51e-07}$ $(\num{9.40e-07})$ & $\num{6.68e-07}$ $(\num{5.23e-07})$ & $\num{8.77e-07}$ $(\num{7.23e-07})$ & $0.13$\\
    & $\num{7.51e-07}$ $(\num{6.82e-07})$ & $\num{4.75e-07}$ $(\num{5.14e-07})$ & $\num{2.19e-07}$ $(\num{2.42e-07})$ & $\num{2.10e-07}$ $(\num{2.01e-07})$ & $0.40$\\ \hline    
    \multirow{2}{*}{$\overline{{\tilde{\varepsilon}}}^{z, r}_0$} & $\num{1.90e-04}$ $(\num{5.71e-05})$ & $\num{1.73e-04}$ $(\num{1.00e-04})$ & $\num{2.55e-03}$ $(\num{9.44e-04})$ & $\num{3.18e-02}$ $(\num{1.40e-02})$ & $-1.45$\\
    & $\num{7.56e-04}$ $(\num{1.78e-04})$ & $\num{2.62e-04}$ $(\num{3.10e-05})$ & $\num{2.71e-04}$ $(\num{1.01e-04})$ & $\num{8.09e-04}$ $(\num{3.76e-04})$ & $0.03$\\ \hline    
    \multirow{2}{*}{$\overline{{\tilde{\varepsilon}}}^{\gamma, r}_0$} & $\num{9.62e-01}$ $(\num{3.44e-02})$ & $\num{1.10e+00}$ $(\num{1.25e-01})$ & $\num{1.95e+00}$ $(\num{3.82e-01})$ & $\num{1.65e+00}$ $(\num{4.19e-01})$ & $-0.20$\\ 
    & $\num{9.78e-03}$ $(\num{6.69e-04})$ & $\num{1.58e-02}$ $(\num{1.09e-03})$ & $\num{4.47e-02}$ $(\num{1.60e-03})$ & $\num{8.87e-02}$ $(\num{2.47e-03})$ & $-0.64$\\ \hline    
    \multirow{2}{*}{$\overline{\tau}$} & $\num{3.57e+03}$ & $\num{1.08e+04}$ & $\num{4.04e+04}$ & $\num{8.03e+04}$ & \\  
    & $\num{1.76e+03}$ & $\num{5.57e+03}$ & $\num{2.10e+04}$ & $\num{4.14e+04}$ & \\ \hline    
   \end{tabular}}
  \end{center}
\caption{Mean relative MSE values, empirical convergence rates of $\left(Y_0, Z_0, \Gamma_0 \right)$ from LDBSDE and DLDBSDE schemes and their average runtimes in Example~\ref{ex3} for $d=50$ and $N \in \{2, 8, 32, 64\}$. The STD of the relative MSE values at $t_0$ is given in the brackets.}
\label{tab4} 
}  
\end{table}
The DLDBSDE scheme yields slightly better approximations for the benchmark value $Y_0$. However, it significantly outperforms the LDBSDE scheme in approximating the benchmark values $\left(Z_0, \Gamma_0\right)$, achieving higher accuracy with less computational time.

\section{Conclusions}
\label{sec5}
In this work, we study a forward differential deep learning approach to solve high-dimensional nonlinear BSDEs. The algorithm aims at overcoming the limitation of the forward deep learning BSDE scheme~\cite{raissi2024forward} that struggle with providing high-accurate gradient approximations. By transforming the BSDE problem into a differential deep learning problem using Malliavin calculus, we solve a system of BSDEs that requires the estimation of the triple of processes $\left(Y, Z, \Gamma \right)$. This triple represents the solution, its gradient, and the Hessian matrix. Our approach involves discretizing the integrals via the Euler-Maruyama method and parameterizing the unknown discrete solution triple using three DNNs. The networks parameters are optimized by globally minimizing a differential learning loss function, defined as a weighted sum of the dynamics of the discretized BSDE system that incorporate local loss functions. Through various high-dimensional examples, we demonstrated that our proposed scheme achieves improved accuracy and computational efficiency compared to the forward deep learning scheme~\cite{raissi2024forward}. This increased performance underscores the potential of our approach as a tool for option pricing and hedging problems in high dimensions.

\bibliographystyle{siamplain}
\bibliography{bibfile}

\end{document}